\documentclass[11pt]{amsart}
\usepackage[latin1]{inputenc}

\usepackage{amsmath}
\usepackage{amsfonts}
\usepackage{color}
\usepackage{amssymb}
\usepackage{amsthm}
\usepackage{graphicx,epsfig}
\usepackage{amscd}
\DeclareMathOperator{\dist}{dist}
\DeclareMathOperator{\re}{Re}

\DeclareMathOperator{\Area}{Area}

\newcommand{\nil}{\mathrm{Nil}_3}
\newcommand{\sol}{\mathrm{Sol}_3}

\newcommand{\R}{\mathbb{R}}
\newcommand{\h}{\mathbb{H}}

\newcommand{\rmT}{\mathrm{T}}

\newcommand{\rmd}{\mathrm{d}}

\newcommand{\cS}{{\mathcal S}}

\newcommand{\cE}{{\mathcal E}}

\newcommand{\cH}{{\mathcal H}}
\newcommand{\cA}{{\mathcal A}}

\newcommand{\cC}{{\mathcal C}}

\newcommand{\cP}{{\mathcal P}}

\newcommand{\cG}{{\mathcal G}}

\begin{document}

\newtheorem{thm}{Theorem}[section]
\newtheorem*{thmintro}{Theorem}
\newtheorem{cor}[thm]{Corollary}
\newtheorem{prop}[thm]{Proposition}
\newtheorem{app}[thm]{Application}
\newtheorem{lemma}[thm]{Lemma}
\newtheorem{notation}[thm]{Notations}
\newtheorem{hypothesis}[thm]{Hypothesis}

\newtheorem{defin}[thm]{Definition}
\newenvironment{defn}{\begin{defin} \rm}{\end{defin}}
\newtheorem{remk}[thm]{Remark}
\newenvironment{rem}{\begin{remk} \rm}{\end{remk}}
\newtheorem{exa}[thm]{Example}
\newenvironment{ex}{\begin{exa} \rm}{\end{exa}}
\newtheorem{cla}[thm]{Claim}
\newenvironment{claim}{\begin{cla} \rm}{\end{cla}}

\title[Half-space theorems]{Half-space theorems for minimal surfaces \\ in $\nil$ and $\sol$}
\author{Beno\^\i t Daniel, William H. Meeks III and Harold Rosenberg}
\date{}


\address{Universit\'e Paris-Est Cr\'eteil, D\'epartement de Math\'ematiques, UFR des Sciences et Technologies, 61 avenue du G\'en\'eral de Gaulle, 94010 Cr\'eteil cedex, FRANCE}
\email{daniel@univ-paris12.fr}
\address{Mathematics Department, University of Massachusetts, Amherst, MA 01003}
\email{profmeeks@gmail.com}
\address{Instituto Nacional de Matem\'atica Pura e Aplicada (IMPA), Estrada Dona Castorina 110, 22460-320 Rio de Janeiro - RJ, BRAZIL}
\email{rosen@impa.br}

\begin{abstract}
We prove some half-space theorems for minimal surfaces in the Heisenberg group $\nil$ and the Lie group $\sol$ endowed with their left-invariant Riemannian metrics. If $\cS$ is a properly immersed minimal surface in $\nil$ that lies on one side of some entire minimal graph $\cG$, then $\cS$ is the image of $\cG$ by a vertical translation. If $\cS$ is a properly immersed minimal surface in $\sol$ that lies on one side of a  special plane  $\cE^t$ (see the discussion just before Theorem~\ref{sol:halfspace} for the definition of a special plane in $\sol$), then $\cS$ is the special plane $\cE^u$ for some $u\in\R$. 
\end{abstract}

\maketitle

\section{Introduction} \label{sec:intro}

A classical theorem in the global theory of proper minimal surfaces in Euclidean $3$-space is the \emph{half-space theorem} by Hoffman and Meeks \cite{HM}: if $\cS$ is a properly immersed minimal surface in $\R^3$ that lies on one side of some plane $\cP$, then $\cS$ is a plane parallel to $\cP$. The proof uses the maximum principle and the fact that catenoids converge to a double cover of a punctured plane as the necksize goes to zero. As a consequence, they proved the \emph{strong half-space theorem}: two properly immersed minimal surfaces in $\R^3$ that do not intersect must be parallel planes.

These theorems have been generalized to some other ambient simply connected homogeneous manifolds. Let us first observe that there is no half-space theorem in Euclidean spaces of dimensions $n\geqslant4$, since there exist rotational proper minimal hypersurfaces contained in a slab.

In hyperbolic $3$-space $\h^3$, one does not have a half-space theorem for minimal surfaces (indeed, for instance any smooth closed curve in the asymptotic boundary of $\h^3$ bounds a minimal surface), but one has half-space theorems for constant mean curvature (CMC) $1$ surfaces \cite{roro}, which can be obtained using rotational catenoid cousins (see also \cite{chr}).
One of the reasons that halfspace theorems exist for CMC 1 surfaces in $\h^3$
is that the ``critical" value for mean curvature in $\h^3$ is $1$, i.e., there exist compact CMC $H$ surfaces in $\h^3$ if and only if $|H|>1$.

Similarly, there is no half-space theorem for minimal surfaces in $\h^2\times\R$, since catenoids (i.e., rotational minimal surfaces) are contained in a slab \cite{NeRo1,NeRo1erratum}. On the other hand, Hauswirth, Rosenberg and Spruck proved a half-space theorem for CMC $\frac12$ surfaces in $\h^2\times\R$.

\begin{thm}[\cite{hrs}]
Let $\cS$ be a properly immersed CMC $\frac12$ surface in $\h^2\times\R$.
\begin{itemize}
\item If $\cS$ is contained on the mean convex side of a horocylinder $\cC$, then $\cS$ is a horocylinder parallel to $\cC$.
\item If $\cS$ is embedded and contains a horocylinder $\cC$ on its mean convex side, then $\cS$ is a horocylinder parallel to $\cC$.
\end{itemize}
(A horocylinder is a product $\gamma\times\R$ where $\gamma\subset\h^2$ is a horocycle.)
\end{thm}

Since rotational CMC $\frac12$ surfaces are not suitable to obtain this theorem, their proof uses a continuous family of compact annuli bounded by two circles in parallel horocylinders, one circle being fixed and the other one having a radius going to infinity. To do this, they use the Schauder fixed point theorem and elliptic PDE techniques.

The aim of this paper is to prove using geometric arguments some half-space theorems for minimal surfaces in two simply connected homogenous $3$-manifolds, the Heisenberg group $\nil$ and the Lie group $\sol$, which are two manifolds admitting isometries with remarkable properties.

The $3$-dimensional Heisenberg group $\nil$ admits a Riemannian submersion $\pi \colon \nil\to\R^2$. Translations along the fibers are isometries called \emph{vertical translations}.

The inverse image by $\pi$ of a straight line in $\R^2$ is a minimal surface called \emph{vertical plane}; two vertical planes are said to be \emph{parallel} if their images by $\pi$ are parallel straight lines. These vertical planes are minimal, stable, isometric to $\R^2$ but not totally geodesic (in fact, there are no local totally geodesic surfaces in $\nil$).

Other examples of stable minimal surfaces in $\nil$ are \emph{entire minimal graphs}, i.e., minimal surfaces $\cG$ such that $\pi_{|\cG} \colon \cG\to\R^2$ is a diffeomorphism. There exist many entire minimal graphs and they were classified by Fernandez and Mira \cite{fmtams}. Some examples (in the usual coordinates $(x_1,x_2,x_3)$ described in Section \ref{sec:nil}) are given by $x_3=0$ (which is rotational) and $x_3=\frac{x_1x_2}2$ (which is invariant by a one-parameter family of translations). In this space, there exist entire minimal graphs of parabolic conformal type and of hyperbolic conformal type.

Using rotational catenoids, Abresch and Rosenberg obtained the following result.

\begin{thm}[\cite{ar2}] \label{thmar}
Let $\cS$ be a properly immersed minimal surface in $\nil$. If $\cS$ lies on one side of the surface of equation $x_3=0$, then $\cS$ is the surface of equation $x_3=c$ for some $c\in\R$.
\end{thm}

Let us observe that one has the analogous statement for the surfaces of equation $x_3=ax_1+bx_2+c$ for any $a$, $b$, $c$ since all these surfaces are congruent.

Daniel and Hauswirth proved a \emph{vertical half-space theorem} for minimal surfaces in $\nil$.

\begin{thm}[\cite{dh}] \label{thmdh}
Let $\cS$ be a properly immersed minimal surface in $\nil$. If $\cS$ lies on one side of some vertical plane $\cP$, then $\cS$ is a vertical plane parallel to $\cP$.
\end{thm}

To do this, they constructed a one-parameter family of \emph{horizontal catenoids}, which are non-rotational properly embedded minimal annuli that converge to a double cover of a punctured vertical plane as the necksize goes to zero (they are semi-explicit and obtained by integrating a Weierstrass-type representation).

In this paper, we will give another proof of this theorem (requiring less computations) and we will prove the following result.

\begin{thm} \label{nil:halfspace}
Let $\cS$ be a properly immersed minimal surface in $\nil$. If $\cS$ lies on one side of some entire minimal graph $\cG$, then $\cS$ is the image of $\cG$ by a vertical translation.
\end{thm}

Theorem \ref{thmar} is the particular case of Theorem \ref{nil:halfspace} when $\cG$ is the surface of equation $x_3=0$.
It is natural to conjecture that two properly immersed minimal surfaces in $\nil$ that do not intersect are either two parallel vertical planes or an entire minimal graph and its image by a vertical translation (this would be the analogue of the strong half-space theorem of $\R^3$).

The Lie group $\sol$ admits a Riemannian submersion $\psi \colon \sol\to\R$ such that, for any $s\in\R$, the surface
$$\cE^s:=\psi^{-1}(s)$$ is minimal, stable, isometric to $\R^2$ but not totally geodesic. We will call these surfaces \emph{special planes}. Other remarkable minimal foliations in $\sol$ are the foliations $(\cH_j^t)_{t\in\mathbb{R}}$ ($j=1,2$) where $\cH_j^t$ is defined by $x_j=t$ in the usual coordinates $(x_1,x_2,s)$ defined in Section \ref{sec:sol}. The leaves are totally geodesic, stable and have intrinsic curvature $-1$ (in fact they are the only totally geodesic surfaces in $\sol$). They are symmetry planes in $\sol$, which permits Alexandrov reflection. Two surfaces $\cH_j^t$ and $\cH_k^u$ are congruent.

In this paper, we will prove the following theorem.

\begin{thm} \label{sol:halfspace}
Let $\cS$ be a properly immersed minimal surface in $\sol$. If $\cS$ lies on one side of some special plane $\cE^t$ ($t\in\R$), then $\cS$ is the special plane $\cE^u$ for some $u\in\R$.
\end{thm}

\begin{rem}
There is no half-space theorem for the surfaces $\cH_1^t$ and $\cH_2^t$; indeed the equation $x_1=ae^{-s}$ for $a\neq0$ defines a properly embedded minimal surface lying on one side of $\cH_1^0$.
\end{rem}

All these half-space theorems are maximum principle at infinity theorems. Let us mention that a maximum principle at infinity was proved by Rosenberg for CMC surfaces with \emph{large} mean curvature in homogeneously regular $3$-manifolds.

\begin{thm}[\cite{rosen}]
Let $N$ be an orientable homogeneously regular $3$-manifold. There exists a constant $c>0$ such that, whenever $H\geqslant c$ and $\cS_1$ and $\cS_2$ are properly embedded CMC $H$ surfaces in $N$ which bound a connected domain $W$, then the mean curvature vector points out of $W$ along the boundary of $W$.
\end{thm}

We now outline the proofs of the theorems. We assume that there exist a properly immersed minimal surface $\cS$ lying on one side of the surface $\Sigma$ with respect to which we want to prove the half-space theorem (i.e., $\Sigma$ is an entire minimal graph in $\nil$, a vertical plane in $\nil$ or a special plane in $\sol$). We assume that $\cS$ and $\Sigma$ are not congruent. We consider the image $\Sigma^\varepsilon$ of $\Sigma$ by a ``small" translation, a fixed circle (``small") in $\Sigma^\varepsilon$ and a circle of varying radius (``large") in $\Sigma$. We construct a least area annulus bounded by the two circles (using the Douglas criterion) and prove using curvature estimates that a subsequence of these annuli converges to a properly embedded annulus bounded by the small circle as the radius of the large circle goes to infinity.

Then we prove that the distance between the limit annulus and the surface $\Sigma$ is positive. To do this we distinguish two cases.
\begin{enumerate}
\item When $\Sigma$ is a vertical plane in $\nil$ or a special plane in $\sol$ (Theorems \ref{thmdh} and \ref{sol:halfspace}), we prove, using curvature estimates, that the limit annulus is quasi-isometric to $\Sigma$ and hence parabolic, and we find a suitable bounded subharmonic function on the limit annulus.
\item When $\Sigma$ is an entire minimal graph in $\nil$ (Theorem \ref{nil:halfspace}), we cannot hope using such an argument since some entire minimal graphs are hyperbolic. Instead, we use a nodal domain argument to prove that the limit annulus is a graph, and then we use a generalization by Leandro and Rosenberg \cite{lr} of a theorem by Collin and Krust \cite{ck} about graphs with prescribed mean curvature. This theorem states that if $u$ and $v$ are two solutions to the same prescribed mean curvature graph equation over a domain $\Omega\subset\R^2$ that coincide on $\partial\Omega$, then either $u-v$ is unbounded or $u\equiv v$ on $\Omega$.
\end{enumerate}

This implies that some annulus bounded by two circles must intersect $\cS$, and we conclude by translating this annulus until reaching an interior last point of contact, contradicting the maximum principle.

The paper is organized as follows. Section \ref{sec:solhalfspace} is devoted to preliminaries about $\sol$ and to the proof of Theorem \ref{sol:halfspace}. Section \ref{sec:nilhalfspace} contains preliminaries about $\nil$ and the proofs of Theorems \ref{thmdh} and \ref{nil:halfspace}. Finally, in the Appendix we include a short proof of a result about subharmonic maps needed in the proofs.

\section{A half-space theorem in $\sol$} \label{sec:solhalfspace}

\subsection{The Lie group $\sol$} \label{sec:sol}

The Lie group $\sol$ can be viewed as $\R^3$ endowed with the Riemannian metric
$$e^{2s}\rmd x_1^2+e^{-2s}\rmd x_2^2+\rmd s^2$$ where $(x_1,x_2,s)$ denote the canonical coordinates of $\R^3$.

In these coordinates, the Riemannian submersion is given by $$\psi \colon (x_1,x_2,s)\mapsto s,$$ and so the special plane $\cE^t$ is simply defined by the equation $s=t$.

We consider the left-invariant orthonormal frame $(E_1,E_2,E_3)$ defined by
$$E_1=e^{-s}\frac{\partial}{\partial x_1},\quad
E_2=e^{s}\frac{\partial}{\partial x_2},\quad
E_3=\frac{\partial}{\partial s}.$$ We call it the \emph{canonical frame}.
The expression of the Riemannian connection $\widehat\nabla$ of $\sol$ in this frame is the following:
\begin{equation} \label{sol:connection}
\begin{array}{lll}
\widehat\nabla_{E_1}E_1=-E_3, &
\widehat\nabla_{E_2}E_1=0, &
\widehat\nabla_{E_3}E_1=0, \\
\widehat\nabla_{E_1}E_2=0, &
\widehat\nabla_{E_2}E_2=E_3, &
\widehat\nabla_{E_3}E_2=0, \\
\widehat\nabla_{E_1}E_3=E_1, &
\widehat\nabla_{E_2}E_3=-E_2, &
\widehat\nabla_{E_3}E_3=0.
\end{array}
\end{equation}

The isometry group of $\sol$ has dimension $3$. The connected component of the identity is generated by the following three families of isometries:
$$(x_1,x_2,s)\mapsto(x_1+c,x_2,s),\quad
(x_1,x_2,s)\mapsto(x_1,x_2+c,s),$$
$$(x_1,x_2,s)\mapsto(e^{-c}x_1,e^cx_2,s+c).$$
The corresponding Killing fields are
$$F_1=\frac\partial{\partial x_1},\quad
F_2=\frac\partial{\partial x_2},\quad
F_3=-x_1\frac\partial{\partial x_1}+x_2\frac\partial{\partial x_2}+\frac\partial{\partial s}.$$ We will call \emph{translations} isometries belonging to the identity component of the identity (they are in fact left multiplications for the Lie group structure).

The isotropy group of the origin $(0,0,0)$ is isomorphic to the dihedral group $\mathrm{D}_4$ and is generated by the following two orientation-reversing isometries:
\begin{equation}
\sigma \colon (x_1,x_2,s)\mapsto(x_2,-x_1,-s),\quad
\tau \colon (x_1,x_2,s)\mapsto(-x_1,x_2,s).
\end{equation}
The reflection with respect to the surface $x_2=0$ is given by $\sigma^2\tau$.

For more details, we refer to \cite{dm} and references therein.

%

\subsection{Proof of Theorem \ref{sol:halfspace}} \label{sec:proofsol}

Before proving the theorem, we will need some preliminary results.

\begin{lemma} \label{sol:subharmonic}
Let $\Sigma$ be a minimal surface (possibly with boundary) in $\sol$ such that $$0<s\leqslant 2$$ on $\Sigma$. Then the function $$\varphi:=\frac1{s}$$ is subharmonic on $\Sigma$.
\end{lemma}

\begin{proof}
We view $\Sigma$ as a conformal minimal immersion $X=(x_1,x_2,s) \colon \Sigma\to\sol$ from a Riemann surface $\Sigma$. Let $z$ be a conformal coordinate. We set $$A_1=e^{-s}x_{1z},\quad A_2=e^{s}x_{2z},\quad A_3=x_{3z},$$ so that
$$X_z=A_1E_1+A_2E_2+A_3E_3.$$ The conformality of $X$ means that
\begin{equation} \label{sol:conformality}
A_1^2+A_2^2+A_3^2=0.
\end{equation}
Since $X$ is minimal, we have $\widehat\nabla_{X_{\bar z}}X_z=0$, and so
\begin{equation} \label{sol:diffA}
A_{3\bar z}=\langle E_3,X_z\rangle_{\bar z}
=\langle\widehat\nabla_{X_{\bar z}}E_3,X_z\rangle
=|A_1|^2-|A_2|^2.
\end{equation}
by \eqref{sol:connection}. Then by \eqref{sol:diffA} we have $$\varphi_{z\bar z}=-\left(\frac{A_3}{s^2}\right)_{\bar z}
=\frac{2|A_3|^2-s(|A_1|^2-|A_2|^2)}{s^3}.$$
On the other hand, by \eqref{sol:conformality} we have
\begin{eqnarray*}
|A_3|^4 & = & |A_1^2+A_2^2|^2=|A_1|^4+|A_2|^4+2\re(A_1^2\bar A_2^2) \\
& \geqslant & |A_1|^4+|A_2|^4-2|A_1|^2|A_2|^2=(|A_1|^2-|A_2|^2)^2.
\end{eqnarray*}
Consequently we get $-|A_3|^2\leqslant|A_1|^2-|A_2|^2\leqslant|A_3|^2$, and so, since $0<s\leqslant2$, we conclude that $\varphi_{z\bar z}\geqslant0$.
\end{proof}

\begin{lemma} \label{sol:diffx3}
There exist positive constants $a$, $b$ and $d$ such that, for all stable minimal surfaces $\Sigma$ (possibly with boundary), for all $p\in\Sigma$ such that $\dist(p,\partial\Sigma)>d$, if $|\langle N(p),E_3\rangle|\leqslant a$ where $N(p)$ denotes a unit normal vector to $\Sigma$ at $p$, then there exist points $q_1$ and $q_2$ in $\Sigma$ such that $s(q_1)-s(q_2)\geqslant b$.
\end{lemma}

\begin{figure}[htbp]
\begin{center}
\input{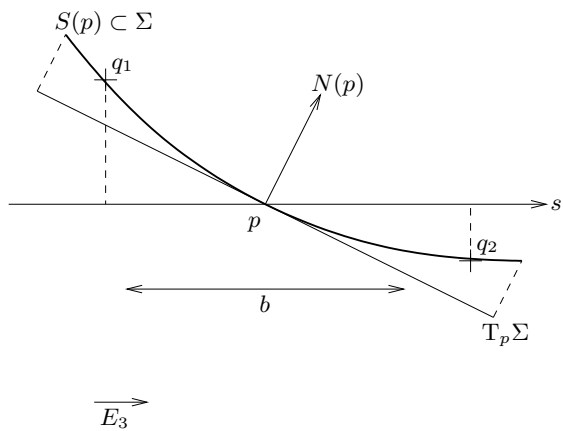}
\caption{A piece of a stable minimal surface.}
\label{figuregraph}
\end{center}
\end{figure}

\begin{proof}
We recall that stable minimal surfaces admit uniform curvature estimates away from their boundary, and so there exist positive constants $d$ and $\delta$ such that, for any a stable minimal surface $\Sigma$ (possibly with boundary), for each $p\in\Sigma$ such that $\dist(p,\partial\Sigma)>d$, there is a piece $S(p)$ of $\Sigma$ around $p$ that is a graph (in exponential coordinates) over the disk in $\rmT_p\Sigma$ of radius $2\delta$ centered at the origin of $\rmT_p\Sigma$. Moreover these graphs have uniformly bounded second fundamental form. If $|\langle N(p),E_3\rangle|$ is smaller than some constant $a>0$, such a piece $S(p)$ necessarily has points $q_1$ and $q_2$ such that $s(q_1)-s(p)\geqslant\frac b2$ and  $s(p)-s(q_2)\geqslant\frac b2$ for some constant $b>0$, and these constants $a$ and $b$ are independent of $\Sigma$ (see Figure \ref{figuregraph}): otherwise one could produce a sequence of such pieces with unbounded second fundamental form. Moreover, these constants are also independent from the point $p$ since $E_3$ and the differences of the $s$ function are invariant by translations.
\end{proof}

From now on, $\cS$ denotes a properly immersed minimal surface in $\sol$. We assume that $\cS$ lies one one side of some special plane $\cE^t$ ($t\in\R$).

Up to isometries in $\sol$, we can assume $t=0$, $s\geqslant 1$ on $\cS$ and
$$\inf\{u\in\R \mid \cS\cap\cE^u\neq\emptyset\}=1.$$

If $\cS\cap\cE^1\neq\emptyset$, then the maximum principle implies that $\cS=\cE^1$. So from now on we assume that
$$\cS\cap\cE^1=\emptyset.$$

For $R>0$ and $h>0$ we let
$$K_R:=\{(x_1,x_2,s)\in\sol \mid x_1^2+x_2^2\leqslant R^2\},$$
$$C_R:=\{(x_1,x_2,s)\in\sol \mid x_1^2+x_2^2=R^2\},$$
$$D_R^h:=K_R\cap\cE^h,$$
$$\Gamma_R^h:=C_R\cap\cE^h,$$
$$Q_R^h:=\bigcup_{t\in[1,1+h]}D_R^t,$$
$$q_R^h:=\bigcup_{t\in[1,1+h]}\Gamma_R^t.$$

We now fix real numbers $r>0$ and $\varepsilon>0$ so that $\varepsilon<1$, $\varepsilon<\frac b2$ where $b$ is defined in Lemma \ref{sol:diffx3},
\begin{equation} \label{areasol}
\Area(q_r^{\varepsilon})<\Area(D_r^1)+\Area(D_r^{1+\varepsilon}).
\end{equation}
and
\begin{equation} \label{sol:intersect}
\cS\cap Q_{er}^{\varepsilon}=\emptyset,
\end{equation}
which is possible since $\cS$ is proper in $\nil$.


\begin{claim} \label{sol:compactannulus}
If $R>r$, then there exists a least area annulus $\cA_R$ bounded by $\Gamma^{1+\varepsilon}_r$ and $\Gamma^1_R$. Moreover this annulus lies between the special planes $\cE^{1+\varepsilon}$ and $\cE^1$ and it is embedded.
\end{claim}

\begin{proof}
The solutions to the Plateau problem for $\Gamma^{1+\varepsilon}_r$ and $\Gamma^1_R$ are respectively $D^{1+\varepsilon}_r$ and $D^1_R$; indeed the unique compact minimal surface bounded by an embedded closed curve in the special plane $\cE^h$ ($h\in\R$) is the part of this special plane bounded by this curve, by the maximum principle (since we have the minimal foliation $(\cE^u)_{u\in\R}$). The total area of these two disks is $\Area(D^1_R)+\Area(D^{1+\varepsilon}_r)$.

Let
$$M:=(D^1_R\setminus D^1_r)\cup q_r^\varepsilon.$$ By \eqref{areasol}, the area of $M$ is smaller than
$\Area(D^1_R)+\Area(D^{1+\varepsilon}_r)$.

Consequently, the Douglas criterion implies the existence of a least area annulus $\cA_R$ bounded by $\Gamma^{1+\varepsilon}_r$ and $\Gamma^1_R$. This annulus is embedded  by the Geometric Dehn's Lemma in \cite{my}.

\end{proof}

By Lemma \ref{sol:diffx3}, if $\Omega\subset\cA_R$ and $\dist(\Omega,\partial\cA_R)>d$, then $|\langle N,E_3\rangle|>a$ where $N$ denotes the unit normal to $\cA_R$ (otherwise there would exist two points $q_1,q_2\in\cA_R$ such that $s(q_1)-s(q_2)\geqslant b>2\varepsilon$, which contradicts the fact that $\cA_R$ lies between the special planes $\cE^{1+\varepsilon}$ and $\cE^1$). In particular, $\Omega$ is transverse to $E_3$.

We now introduce a constant $\rho>r$ such that
\begin{equation} \label{bounddist}
\dist(\Gamma_r^{1+\varepsilon},q_\rho^\varepsilon)>d.
\end{equation}
If $R$ is large enough, this implies that $\cA_R$ intersects transversely $q_\rho^\varepsilon$ in a smooth curve $\widetilde\Gamma_R$. We denote by $\widetilde\cA_R$ the part of $\cA_R$ lying outside $Q_\rho^\varepsilon$: it is an annulus bounded by $\widetilde\Gamma_R$ and $\Gamma_R^{1+\varepsilon}$.

\begin{claim} \label{sol:limitannulus}
Let $(R_n)$ be an increasing sequence of positive real numbers such that $R_n\to+\infty$ as $n\to+\infty$. Then, up to a subsequence, the annuli $\widetilde\cA_{R_n}$ converge to a properly embedded minimal annulus $\widetilde\cA_\infty$ whose boundary is a closed curve $\widetilde\Gamma_\infty\subset q_\rho^\varepsilon$ and
such that $$\inf_{\widetilde\cA_\infty}s>1.$$
\end{claim}


\begin{proof}
The annuli $\cA_{R_n}$ are stable so they admit uniform curvature bounds for their points lying at a distance $\geqslant d$ from their boundary. Hence, by \eqref{bounddist}, up to a subsequence, the annuli $\widetilde\cA_{R_n}$ converge to a properly embedded minimal surface $\widetilde\cA_\infty$ whose boundary is a curve $\widetilde\Gamma_\infty$ lying between the special planes $\cE^{1+\varepsilon}$ and $\cE^1$ and in $q_\rho^\varepsilon$ (properness of the limit annulus follows easily deduced from the fact that outside a fixed sized regular neighborhood of  its boundary, each $\cA_{R_n}$ is transverse to $E_3$ and graphical over its $s$-projection to $\cE^1$). Moreover, since $\widetilde\cA_{R_m}$ lies above $\widetilde\cA_{R_n}$ if $m>n$, the curve $\widetilde\Gamma_\infty$ lies above all the curves $\widetilde\Gamma_{R_n}$; in particular we have $$h:=\min_{\widetilde\Gamma_\infty}s>1.$$

Let $N$ be a unit normal vector field to $\widetilde\cA_\infty$. Then we have $|\langle N,E_3\rangle|\geqslant a$, so $\widetilde\cA_\infty$ can be written as a graph $s=f(x_1,x_2)$ over $$L:=\{(x_1,x_2,1)\in\cE^1 \mid x_1^2+x_2^2\geqslant\rho^2\},$$
where $f$ is a function with bounded gradient. The vector field $E_3$ is orthogonal to $\cE^1$, and by construction, the annulus $\widetilde\cA_\infty$ lies between the special planes $\cE^1$ and $\cE^{1+\varepsilon}$. This implies that the map $f \colon L\to\widetilde\cA_\infty$ is a quasi-isometry. Since $\cE^1$ is flat, this implies that $\widetilde\cA_\infty$ is parabolic (it has quadratic area growth).

On $\widetilde\cA_\infty$ we have $1\leqslant s\leqslant1+\varepsilon\leqslant2$, so by Lemma \ref{sol:subharmonic} the function $\varphi:=1/s$ is subharmonic on $\widetilde\cA_\infty$. Moreover, $1/(1+\varepsilon)\leqslant\varphi\leqslant1$ on $\widetilde\cA_{\infty}$, so since $\widetilde\cA_\infty$ has parabolic conformal type, Proposition \ref{bounded} implies that
$$\varphi\leqslant\sup_{\partial\widetilde\cA_\infty}\varphi=\frac1h$$ on $\widetilde\cA_\infty$. Consequently we have $s\geqslant h>1$ on $\widetilde\cA_\infty$.
\end{proof}

We can now conclude the proof of the theorem.

\begin{figure}[htbp]
\begin{center}
\input{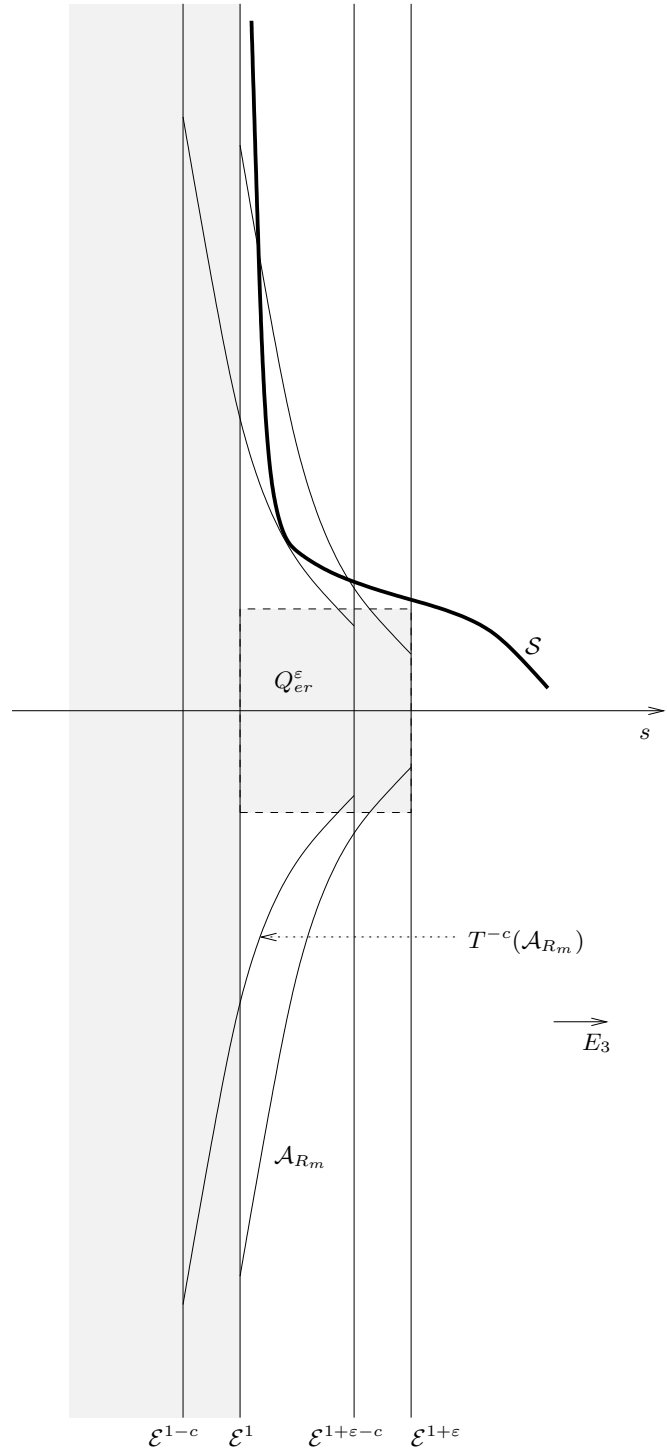}
\caption{Translating the annulus $\cA_{R_m}$.}
\label{figuretranslation}
\end{center}
\end{figure}

\begin{proof}[Proof of Theorem \ref{sol:halfspace}]
Because of Claim \ref{sol:limitannulus}, there exists $m\in\mathbb{N}$ such that $\cS\cap\widetilde\cA_{R_m}\neq\emptyset$, and so
$$\cS\cap\cA_{R_m}\neq\emptyset.$$ 

For $c\in\R$, we let $T^c \colon \sol\to\sol$ denote the isometry $(x_1,x_2,s)\mapsto(e^{-c}x_1,e^cx_2,s+c)$.
We consider the annuli $T^{-c}(\cA_{R_m})$ for $c\geqslant0$. We notice that $\cS\cap T^{-c}(\cA_{R_m})=\emptyset$ when $c>\varepsilon$. Then there exists a largest $c$ for which $$\cS\cap T^{-c}(\cA_{R_m})\neq\emptyset.$$

We claim that no point of intersection lies on the boundary of $T^{-c}(\cA_{R_m})$. Indeed, this boundary consists of $T^{-c}(\Gamma_r^{1+\varepsilon})$ and $T^{-c}(\Gamma_{R_m}^1)$. The curve $T^{-c}(\Gamma_r^{1+\varepsilon})$ is defined by
$$\left\{\begin{array}{rcl}
e^{-2c}x_1^2+e^{2c}x_2^2 & = & r^2, \\
s & = & 1+\varepsilon-c,\end{array}\right.$$
so, since $0\leqslant c\leqslant\varepsilon<1$, it is contained in $Q_{er}^{\varepsilon}$ and hence cannot intersect $\cS$ by \eqref{sol:intersect}, On the other hand, $T^{-c}(\Gamma_{R_m}^1)$ lies below $\cE^1$ and hence cannot intersect $\cS$ either.

Consequently there exists an intersection point of $\cS$ and $T^{-c}(\cA_{R_m})$ lying in the interior of $T^{-c}(\cA_{R_m})$. But since $c$ is maximal, $\cS$ lies on one side of $T^{-c}(\cA_{R_m})$; this contradicts the maximum principle (see Figure \ref{figuretranslation}).
\end{proof}

%
%

\section{Half-space theorems in $\nil$} \label{sec:nilhalfspace}

\subsection{The Lie group $\nil$} \label{sec:nil}

The $3$-dimensional Heisenberg group $\nil$ can be viewed as $\R^3$ endowed with the metric
$$\rmd x_1^2+\rmd x_2^2+\left(\frac12(x_2\rmd x_1-x_1\rmd x_2)+\rmd x_3\right)^2.$$
The projection $\pi \colon \nil\to\R^2,(x_1,x_2,x_3)\mapsto(x_1,x_2)$ is a Riemannian submersion.

We consider the left-invariant orthonormal frame $(E_1,E_2,E_3)$ defined by
$$E_1=\frac\partial{\partial x_1}-\frac{x_2}2\frac\partial{\partial x_3},\quad
E_2=\frac\partial{\partial x_2}+\frac{x_1}2\frac\partial{\partial x_3},\quad
E_3=\xi=\frac\partial{\partial x_3}.$$ We call it the \emph{canonical frame}.
The expression of the Riemannian connection $\widehat\nabla$ of $\nil$ in this frame is the following:
\begin{equation} \label{nil:connection}
\begin{array}{lll}
\widehat\nabla_{E_1}E_1=0, &
\displaystyle{\widehat\nabla_{E_2}E_1=-\frac12E_3}, &
\displaystyle{\widehat\nabla_{E_3}E_1=-\frac12E_2}, \\
\displaystyle{\widehat\nabla_{E_1}E_2=\frac12E_3}, &
\widehat\nabla_{E_2}E_2=0, &
\displaystyle{\widehat\nabla_{E_3}E_2=\frac12E_1}, \\
\displaystyle{\widehat\nabla_{E_1}E_3=-\frac12E_2}, &
\displaystyle{\widehat\nabla_{E_2}E_3=\frac12E_1}, &
\widehat\nabla_{E_3}E_3=0.
\end{array}
\end{equation}

A vector is said to be vertical if it is proportional to $\xi$, and horizontal if it is orthogonal to $\xi$. A surface is said to be a (local) \emph{$\xi$-graph} if it is transverse to $\xi$. We will call the inverse image by $\pi$ of a straight line in $\R^2$ a  \emph{vertical plane}.

The isometry group of $\nil$ has dimension $4$. The connected component of the identity is generated by the following four families of isometries:
$$(x_1,x_2,x_3)\mapsto\left(x_1+c,x_2,x_3+\frac{cx_2}2\right),$$
$$(x_1,x_2,x_3)\mapsto\left(x_1,x_2+c,x_3-\frac{cx_1}2\right),$$
$$(x_1,x_2,x_3)\mapsto(x_1,x_2,x_3+c),$$
$$(x_1,x_2,x_3)\mapsto((\cos\theta)x_1-(\sin\theta)x_2,(\sin\theta)x_1+(\cos\theta)x_2,x_3).$$
The corresponding Killing fields are
$$F_1=\frac\partial{\partial x_1}+\frac{x_2}2\frac\partial{\partial x_3},\quad
F_2=\frac\partial{\partial x_2}-\frac{x_1}2\frac\partial{\partial x_3},\quad
F_3=\xi=\frac\partial{\partial x_3},$$
$$F_4=(-(\sin\theta)x_1-(\cos\theta)x_2)\frac\partial{\partial x_1}
+((\cos\theta)x_1-(\sin\theta)x_2)\frac\partial{\partial x_2}.$$
We will call \emph{translations} isometries generated by the first three of these families (they are in fact left multiplications for the Lie group structure).

The entire isometry group of $\nil$ is generated by the aforementioned isometries and the orientation-reversing isometry
$$(x_1,x_2,x_3)\mapsto(-x_1,x_2,-x_3).$$

For more details, we refer to \cite{dh} and references therein.

\subsection{A new proof of Theorem \ref{thmdh}} \label{sec:proofdh}

We can prove Theorem \ref{thmdh} in a very similar way to that of Theorem \ref{sol:halfspace}. For our purpose it will be useful to introduce the following coordinates in $\nil$:
$$y_1=x_1,\quad y_2=x_2,\quad y_3=x_3+\frac{x_1x_2}2.$$ In these coordinates we have
$$E_1=\frac\partial{\partial y_1},\quad
E_2=\frac\partial{\partial y_2}+y_1\frac\partial{\partial y_3},\quad
E_3=\xi=\frac\partial{\partial y_3}.$$
What is important is that $E_1$ is the partial derivative with respect to a coordinate whose level sets are precisely the vertical planes to which $E_1$ is orthogonal. The proof of Theorem \ref{thmdh} will be analogous to that of Theorem \ref{sol:halfspace}, replacing the $s$ coordinate in $\sol$ by the $y_1$ coordinate in $\nil$. The main difference lies in the analogue of Lemma \ref{sol:subharmonic}, which is the following.

\begin{lemma} \label{nil:subharmonic}
Let $\Sigma$ be a minimal surface (possibly with boundary) in $\nil$ such that $$0<y_1\leqslant 4$$ on $\Sigma$. Then the function $$\varphi:=\frac1{y_1}$$ is subharmonic on $\Sigma$.
\end{lemma}

\begin{proof}
We view $\Sigma$ as a conformal minimal immersion $X=(y_1,y_2,y_3) \colon \Sigma\to\nil$ from a Riemann surface $\Sigma$. Let $z$ be a conformal coordinate. We set
$$A_1=y_{1z},\quad A_2=y_{2z},\quad A_3=y_{3z}-y_1y_{2z},$$ so that
$$X_z=A_1E_1+A_2E_2+A_3E_3.$$ The conformality of $X$ means that
\begin{equation} \label{nil:conformality}
A_1^2+A_2^2+A_3^2=0.
\end{equation}
Since $X$ is minimal, we have $\widehat\nabla_{X_{\bar z}}X_z=0$, and so
\begin{equation} \label{nil:diffA}
A_{1\bar z}=\langle E_1,X_z\rangle_{\bar z}
=\langle\widehat\nabla_{X_{\bar z}}E_1,X_z\rangle
=-\frac12(\bar A_2A_3+A_2\bar A_3)
\end{equation}
by \eqref{nil:connection}. Then by \eqref{nil:diffA} we have
$$\varphi_{z\bar z}=-\left(\frac{A_1}{y_1^2}\right)_{\bar z}
=\frac{4|A_1|^2+y_1(\bar A_2A_3+A_2\bar A_3)}{2y_1^3}.$$
On the other hand, by \eqref{nil:conformality} we have
$$|A_1|^4=|A_2^2+A_3^2|^2=|A_2|^4+|A_3|^4+\bar A_2^2 A_3^2+A_2^2\bar A_3^2$$
and $$(\bar A_2A_3+A_2\bar A_3)^2=\bar A_2^2 A_3^2+A_2^2\bar A_3^2+2|A_2|^2|A_3|^2.$$
Consequently we get $-|A_1|^2\leqslant\bar A_2A_3+A_2\bar A_3\leqslant|A_1|^2$, and so, since $0<y_1\leqslant4$, we conclude that $\varphi_{z\bar z}\geqslant0$.
\end{proof}

\begin{lemma} \label{nil:diffy1}
There exists positive constants $a$, $b$ and $d$ such that, for all stable minimal surface $\Sigma$ (possibly with boundary), for all $p\in\Sigma$ such that $\dist(p,\partial\Sigma)>d$, if $|\langle N(p),E_1\rangle|\leqslant a$ where $N(p)$ denotes a unit normal vector to $\Sigma$ at $p$, then there exist points $q_1$ and $q_2$ in $\Sigma$ such that $y_1(q_1)-y_1(q_2)\geqslant b$.
\end{lemma}

\begin{proof}
The proof is analogous to that of that of Lemma \ref{sol:diffx3}.
\end{proof}

From now on, $\cS$ denotes a properly immersed minimal surface in $\nil$ and $\cP$ a vertical plane. We assume that $\cS$ lies one one side of $\cP$.

For $h\in\R$, we let $\cP^h$ denote the plane of equation $y_1=h$. Up to isometries in $\nil$, we can assume $\cP=\cP^1$, $y_1\geqslant 1$ on $\cS$ and
$$\inf\{u\in\R \mid \cS\cap\cP^u\neq\emptyset\}=1.$$

If $\cS\cap\cP^1\neq\emptyset$, then the maximum principle implies that $\cS=\cP^1$. So from now on we assume that
$$\cS\cap\cP^1=\emptyset.$$

For $R>0$ and $h>0$ we let
$$K_R:=\{(y_1,y_2,y_3)\in\nil \mid y_2^2+y_3^2\leqslant R^2\},$$
$$C_R:=\{(y_1,y_2,y_3)\in\nil \mid y_2^2+y_3^2=R^2\},$$
$$D_R^h:=K_R\cap\cP^h,$$
$$\Gamma_R^h:=C_R\cap\cP^h,$$
$$Q_R^h:=\bigcup_{t\in[1,1+h]}D_R^t,$$
$$q_R^h:=\bigcup_{t\in[1,1+h]}\Gamma_R^t.$$

We now fix real numbers $r>0$ and $\varepsilon>0$ so that $\varepsilon<1$, $\varepsilon<\frac b2$ where $b$ is defined in Lemma \ref{nil:diffy1},
\begin{equation} \label{areanildh}
\Area(q_r^{\varepsilon})<\Area(D_r^1)+\Area(D_r^{1+\varepsilon}).
\end{equation}
and
\begin{equation} \label{nildh:intersect}
\cS\cap Q_{2r}^{\varepsilon}=\emptyset,
\end{equation}
which is possible since $\cS$ is proper in $\nil$.

The following facts can be proved in the same way as Claims \ref{sol:compactannulus} and \ref{sol:limitannulus}.

\begin{claim} \label{nildh:compactannulus}
If $R>r$, then there exists a least area annulus $\cA_R$ bounded by $\Gamma^{1+\varepsilon}_r$ and $\Gamma^1_R$. Moreover this annulus lies between the planes $\cP^{1+\varepsilon}$ and $\cP^1$ and it is embedded.
\end{claim}

By Lemma \ref{nil:diffy1}, if $\Omega\subset\cA_R$ and $\dist(\Omega,\partial\cA_R)>d$, then $|\langle N,E_1\rangle|>a$ where $N$ denotes the unit normal to $\cA_R$ (otherwise there would exist two points $q_1,q_2\in\cA_R$ such that $y_1(q_1)-y_1(q_2)\geqslant b>2\varepsilon$, which contradicts the fact that $\cA_R$ lies between the planes $\cP^{1+\varepsilon}$ and $\cP^1$). In particular, $\Omega$ is transverse to $E_1$.

We now introduce a constant $\rho>r$ such that
\begin{equation} \label{nil:bounddist}
\dist(\Gamma_r^{1+\varepsilon},q_\rho^\varepsilon)>d.
\end{equation}
If $R$ is large enough, this implies that $\cA_R$ intersects transversely $q_\rho^\varepsilon$ in a smooth curve $\widetilde\Gamma_R$. We denote by $\widetilde\cA_R$ the part of $\cA_R$ lying outside $Q_\rho^\varepsilon$: it is an annulus bounded by $\widetilde\Gamma_R$ and $\Gamma_R^{1+\varepsilon}$.

\begin{claim} \label{nildh:limitannulus}
Let $(R_n)$ be an increasing sequence of positive real numbers such that $R_n\to+\infty$ as $n\to+\infty$. Then, up to a subsequence, the annuli $\widetilde\cA_{R_n}$ converge to a properly embedded minimal annulus $\widetilde\cA_\infty$ whose boundary is a closed curve $\widetilde\Gamma_\infty\subset q_\rho^\varepsilon$ and
such that $$\inf_{\widetilde\cA_\infty}y_1>1.$$
\end{claim}

We can now conclude the proof of the theorem.

\begin{proof}[Proof of Theorem \ref{thmdh}]
This is similar to the conclusion to the proof of Theorem \ref{sol:halfspace}. There exists $m\in\mathbb{N}$ such that
$$\cS\cap\cA_{R_m}\neq\emptyset.$$
For $c\in\R$, we let $\Phi^c \colon \nil\to\nil$ denote the isometry $(y_1,y_2,y_3)\mapsto(y_1+c,y_2,y_3+cy_2)$.
We consider the largest $c\geqslant0$ such that $\cS\cap\Phi^{-c}(\cA_{R_m})\neq\emptyset$. We have $c\leqslant\varepsilon<1$, so using \eqref{nildh:intersect} we prove that no intersection points lies on the boundary, and obtain a contradiction with the maximum principle.
\end{proof}

\subsection{Proof of Theorem \ref{nil:halfspace}} \label{sec:proofnil}

In this section, $\cS$ denotes a properly immersed minimal surface in $\sol$ and $\cG$ an entire minimal $\xi$-graph. We assume that $\cS$ lies one one side of $\cG$. For $c\in\R$, we let $T^c \colon \nil\to\nil$ denote the vertical translation $(x_1,x_2,x_3)\mapsto(x_1,x_2,x_3+c)$.

Without loss of generality, we can assume that $\cS$ lies above $\cG$. Let
$$c_0:=\inf\{c\in\R \mid \cS\cap T^c(\cG)\neq\emptyset\}.$$ If $\cS\cap T^{c_0}(\cG)\neq\emptyset$, then the maximum principle implies that $\cS=T^{c_0}(\cG)$. So from now on we assume that $$\cS\cap T^{c_0}(\cG)=\emptyset.$$ Moreover, applying a vertical translation to $\cS$, we can assume that $$c_0=0.$$

For $R>0$ and $h>0$ we let
$$K_R:=\{(x_1,x_2,x_3)\in\nil \mid  x_1^2+x_2^2\leqslant R^2\},$$
$$C_R:=\{(x_1,x_2,x_3)\in\nil \mid  x_1^2+x_2^2=R^2\},$$
$$D_R^h:=K_R\cap T^h(\cG),$$
$$\Gamma_R^h:=C_R\cap T^h(\cG),$$
$$Q_R^h:=\bigcup_{t\in[0,h]}D_R^t,$$
$$q_R^h:=\bigcup_{t\in[0,h]}\Gamma_R^t.$$

We now fix a real number $r\in(0,1)$. The part of $\cG$ bounded by $\Gamma_r^0$ and $\Gamma_1^0$ is stable, hence a small perturbation of its boundary also bounds a stable minimal surface that is near it; in particular, for a small $\varepsilon>0$, there exists a stable annulus $\cA$ bounded by $\Gamma^{\varepsilon}_r$ and $\Gamma^0_1$, and this annulus is a $\xi$-graph. We fix such an $\varepsilon$ and we moreover assume that $\varepsilon$ is small enough so that
\begin{equation} \label{areanil}
\Area(q_r^\varepsilon)<2\Area(D_r^0).
\end{equation}
and
\begin{equation} \label{nil:intersect}
\cS\cap Q_1^{\varepsilon}=\emptyset,
\end{equation}
which is possible since $\cS$ is proper in $\nil$.

\begin{claim} \label{nil:compactannulus}
If $R>r$, then there exists a least area annulus $\cA_R$ bounded by $\Gamma^{\varepsilon}_r$ and $\Gamma^0_R$. Moreover this annulus lies between the graphs $T^\varepsilon(\cG)$ and $\cG$ and it is embedded.
\end{claim}

\begin{proof}
The solutions to the Plateau problem for $\Gamma^{\varepsilon}_r$ and $\Gamma^0_R$ are respectively $D^{\varepsilon}_r$ and $D^0_R$; indeed the unique compact minimal surface bounded by an embedded closed curve in the minimal entire graph $T^h(\cG)$ ($h\in\R$) is the part of this graph bounded by this curve, by the maximum principle (since we have the minimal foliation $(T^c(\cG))_{c\in\R}$). The total area of these two disks is $\Area(D^0_R)+\Area(D^0_r)$.

Let 
$$M:=(D^0_R\setminus D^0_r)\cup q_r^\varepsilon.$$ By \eqref{areanil}, the area of $M$ is smaller than $\Area(D^0_R)+\Area(D^0_r)$.

Consequently, the Douglas criterion implies the existence of a least area annulus $\cA_R$ bounded by $\Gamma^{1+\varepsilon}_r$ and $\Gamma^1_R$. This annulus is embedded  by the Geometric Dehn's Lemma in \cite{my}.

\end{proof}

For $R>r$, we set
$$U_R:=\{(x_1,x_2)\in\R^2 \mid r^2<x_1^2+x_2^2<R^2\}.$$ We also set
$$U_\infty:=\{(x_1,x_2)\in\R^2 \mid r^2<x_1^2+x_2^2\}.$$

\begin{claim} \label{nil:annulusgraph}
If $R>1$, then the annulus $\cA_R$ is a $\xi$-graph of a function $u_R:\overline{U_R}\to\R$.
\end{claim}

\begin{proof}
We first prove that $\cA_R$ is a $\xi$-graph near $\Gamma^{\varepsilon}_r$. The annulus $\cA_R$ lies below $T^\varepsilon(\cG)$; moreover, since $R>1$, $\cA_R$ is situated above $\cA$, which is a $\xi$-graph. As $\cA_R$ is analytic and without branch points up to the boundary \cite{my}, this implies that $\cA_R$ is a $\xi$-graph near $\Gamma^{\varepsilon}_r$ (see Figure \ref{figureannuli}). We now let $N$ be the unit normal vector field to $\cA_R$ so that $N$ points upwards near $\Gamma^\varepsilon_r$.

\begin{figure}[htbp]
\begin{center}
\input{annuli.pstex_t}
\caption{The annulus $\cA_R$.}
\label{figureannuli}
\end{center}
\end{figure}

We now prove that $\cA_R$ is a $\xi$-graph near $\Gamma^0_R$. This comes from the fact that the mean curvature vector of $C_R$ does not vanish and points inside $K_R$, so $\cA_R$ lies inside $K_R$ and cannot be tangent to $C_R$ along $\Gamma^0_R$. Moreover, since $\cA_R$ is embedded, this implies that $N$ also points upwards near $\Gamma^0_R$. 

We now prove that the whole $\cA_R$ is a $\xi$-graph. Assume this is not the case. Then the Jacobi function $\nu:=\langle N,\xi\rangle$ admits a nodal domain $\Omega$ on which $\nu<0$. Denoting by $\lambda_1$ the first eigenvalue of the Jacobi operator, this means that $\lambda_1(\bar\Omega)=0$. On the other hand, $\bar\Omega$ is contained in the interior of $\cA_R$ because $\nu>0$ near the boundary of $\cA_R$. From this we conclude that $$\lambda_1(\cA_R)<\lambda_1(\bar\Omega)=0,$$ which contradicts the fact that $\cA_R$ is stable.
\end{proof}

\begin{claim} \label{nil:limitannulus}
Let $(R_n)$ be an increasing sequence of positive real numbers such that $R_n\to+\infty$ as $n\to+\infty$. Then, up to a subsequence, the functions $u_{R_n}$ converge (in the $\mathrm{C}^2$ topology on compact sets) to a smooth function $u_\infty \colon U_\infty\to\R$. Moreover, this function $u_\infty$ extends to a continuous function
$$u_\infty \colon \overline{U_\infty}\to\R$$ by setting
$$\forall p\in\partial U_\infty,\quad u_\infty(p):=v(p),$$
where $v \colon \R^2\to\R$ is the function whose $\xi$-graph is $T^\varepsilon(\cG)$.
\end{claim}

\begin{proof}
The annuli $\cA_{R_n}$ are stable so one has uniform curvature bounds over any compact subset of $U_\infty$. Hence, up to a subsequence, the annuli $\mathring\cA_{R_n}:=\cA_{R_n}\setminus\partial\cA_{R_n}$ converge to a properly embedded open minimal surface $\mathring\cA_\infty$ lying between the entire graphs $T^\varepsilon(\cG)$ and $\cG$.

Let $N$ be a unit normal vector field to $\mathring\cA_\infty$ and $\nu:=\langle N,\xi\rangle$. Since $\mathring\cA_\infty$ is the limit of the annuli $\mathring\cA_{R_n}$, which are $\xi$-graphs by Claim \ref{nil:annulusgraph}, we have either $\nu\geqslant0$ or $\nu\leqslant0$. Up to a change of orientation, we can assume that $\nu\geqslant0$.
If $\nu$ vanishes at some interior point, then, since $\nu\geqslant0$ and $\nu$ satisfies an elliptic equation of the form $\Delta\nu+V\nu=0$ for some potential $V$, the maximum principle implies that $\nu\equiv0$. This means that $\mathring\cA_\infty$ is part of a vertical surface, hence a vertical plane since $\mathring\cA_\infty$ is minimal; this is a contradiction.

Consequently, $\mathring\cA_\infty$ is the $\xi$-graph of a function $u_\infty \colon U_\infty\to\R$. Finally, the barriers $T^\varepsilon(\cG)$ and $\cA$ imply that $u_\infty$ extends to a continuous function
$$u_\infty \colon \overline{U_\infty}\to\R$$ by setting $u_\infty(p):=v(p)$ for all $p\in\partial U_\infty$.
\end{proof}

\begin{claim} \label{nil:limitannulus2}
We have $$u_\infty\equiv v.$$ In other words, the annulus $\mathring\cA_\infty$ is the part of the entire graph $T^\varepsilon(\cG)$ lying outside $\Gamma^\varepsilon_r$.
\end{claim}

\begin{proof}
The functions $u_\infty$ and $v$ satisfy the minimal graph equation on $U_\infty$, $u_\infty\equiv v$ on $\partial U_\infty$, and $-\varepsilon<u_\infty-v\leqslant 0$ since $\mathring\cA_\infty$ lies between the graphs $T^\varepsilon(\cG)$ and $\cG$. Then Theorem 5.1 in \cite{lr} implies that $u_\infty\equiv v$ (observe that in the theorems of \cite{ck,lr} the functions only need to be continuous along the boundary of the domain).
\end{proof}

We can now conclude the proof of the theorem.

\begin{proof}[Proof of Theorem \ref{nil:halfspace}]
Since the annuli $\cA_{R_n}$ converge, as $n\to+\infty$, to the part of $T^\varepsilon(\cG)$ lying outside $\Gamma^\varepsilon_r$, there exists $m\in\mathbb{N}$ such that
$$\cS\cap\cA_{R_m}\neq\emptyset.$$

We consider the annuli $T^{-c}(\cA_{R_m})$ for $c\geqslant0$. We notice that $\cS\cap T^{-c}(\cA_{R_m})=\emptyset$ when $c>\varepsilon$. Then there exists a largest $c$ for which $$\cS\cap T^{-c}(\cA_{R_m})\neq\emptyset.$$

We claim that no point of intersection lies on the boundary of $T^{-c}(\cA_{R_m})$. Indeed, this boundary consists of $\Gamma_r^{\varepsilon-c}$, which is contained in $Q_1^{\varepsilon}$ and hence cannot intersect $\cS$ by \eqref{nil:intersect}, and of $\Gamma_{R_m}^{-c}$, which lies below $\cG$ and hence cannot intersect $\cS$ either.

Consequently there exists an intersection point of $\cS$ and $T^{-c}(\cA_{R_m})$ lying in the interior of $T^{-c}(\cA_{R_m})$. But since $c$ is maximal, $\cS$ lies on one side of $T^{-c}(\cA_{R_m})$; this contradicts the maximum principle.
\end{proof}

\begin{rem}
Since entire minimal graphs in $\nil$ and entire CMC $\frac12$ graphs in $\h^2\times\R$ are sister surfaces \cite{dh}, it would be interesting to prove a half-space theorem in $\h^2\times\R$ for CMC $\frac12$ surfaces with respect to an entire CMC $\frac12$ graph. However, in this setting there is no known Collin-Krust type theorem (and actually such a theorem fails for \emph{minimal} graphs in $\h^2\times\R$).
\end{rem}

\section{Appendix}


\begin{prop} \label{bounded}
Let $\Sigma$ be a parabolic surface with non-empty boundary. Let $f \colon \Sigma\to\R$ be a bounded subharmonic function. Then $$f\leqslant\sup_{\partial\Sigma}f.$$
\end{prop}

\begin{proof}
Up to lifting $f$ to the universal covering space of $\Sigma$, we may assume that $\Sigma$ is simply connected and so, since it is parabolic, $\Sigma$ is the radius one closed disc with a measure zero set removed from the unit circle.

Let $p\in\Sigma$. Since there is a conformal map sending $p$ to $0$ and leaving the unit circle invariant, we can assume that $p=0$. Since $f$ is subharmonic, for all $r\in(0,1)$ we have
$$f(p)\leqslant\frac1{2\pi r}\int_{|z|=r}f.$$ Since $f$ is bounded, as $r\to1$, these averages converge to the average of $f$ over the unit circle, which is less than or equal to $\sup_{\partial\Sigma}f$. Hence $f(p)\leqslant\sup_{\partial\Sigma}f$. 
\end{proof}

\bibliographystyle{plain}
\bibliography{dmr}

\end{document}